\documentclass[12pt]{article}%
\usepackage{amsfonts}
\usepackage{sw20bams}
\usepackage{amsmath}
\usepackage{amssymb}
\usepackage{graphicx}
\usepackage{endnotes}%
\providecommand{\U}[1]{\protect\rule{.1in}{.1in}}
\let\footnote=\endnote
\begin{document}

\title{Traffic Light Queues and the Poisson Clumping Heuristic}
\author{Steven Finch and Guy Louchard}
\date{January 21, 2019}
\maketitle

\begin{abstract}
In discrete time, $\ell$-blocks of red lights are separated by $\ell$-blocks
of green lights. Cars arrive at random. \ We seek the distribution of maximum
line length of idle cars, and justify conjectured probabilistic asymptotics
for $2\leq\ell\leq3$.

\end{abstract}

Assorted expressions emerge for a certain traffic light problem \cite{Fi1-heu}%
. \ Let $\ell\geq1$ be an integer. \ Let $X_{0}=0$ and $X_{1}$, $X_{2}$,
\ldots, $X_{n}$ be a sequence of independent random variables satisfying%
\[%
\begin{array}
[c]{ccccc}%
\mathbb{P}\left\{  X_{i}=1\right\}  =p, &  & \mathbb{P}\left\{  X_{i}%
=0\right\}  =q &  & \text{if }i\equiv1,2,\ldots,\ell\,\operatorname{mod}%
\,2\ell;
\end{array}
\]%
\[%
\begin{array}
[c]{ccccc}%
\mathbb{P}\left\{  X_{i}=0\right\}  =p, &  & \mathbb{P}\left\{  X_{i}%
=-1\right\}  =q &  & \text{if }i\equiv\ell+1,\ell+2,\ldots,2\ell
\,\operatorname{mod}\,2\ell
\end{array}
\]
for each $1\leq i\leq n$. \ Define $S_{0}=X_{0}$ and $S_{j}=\max\left\{
S_{j-1}+X_{j},0\right\}  $ for all $1\leq j\leq n$. \ Thus cars arrive at a
one-way intersection according to a Bernoulli($p$) distribution; \ when the
signal is red ($1\leq i\leq\ell$), no cars may leave; \ when the signal is
green ($\ell+1\leq i\leq2\ell$), a car must leave (if there is one). The
quantity $M_{n}=\max\nolimits_{0\leq j\leq n}S_{j}$ is the worst-case traffic
congestion (as opposed to the average-case often cited). \ Only the
circumstance when $\ell=1$ is amenable to rigorous treatment \cite{Fi2-heu},
as far as is known. \ We assume that $p<q$ throughout.

The Poisson clumping heuristic \cite{Ald-heu}, while not a theorem, gives
results identical to exact asymptotic expressions when such exist, and
evidently provides excellent predictions otherwise. \ Consider an irreducible
positive recurrent Markov chain with stationary distribution $\pi$. \ For
sufficiently large $k$, the maximum of the chain satisfies%
\[
\mathbb{P}\left\{  M_{n}<k\right\}  \sim\exp\left(  -\frac{\pi_{k}}%
{\mathbb{E}(C)}n\right)
\]
as $n\rightarrow\infty$, where $C$ is the sojourn time in $k$ during a clump
of nearby visits to $k$.

Here is a simple example: for an asymmetric random walk with weak reflection
at the origin, we have%
\[%
\begin{array}
[c]{ccc}%
\pi_{j}=p\pi_{j-1}+q\pi_{j+1}, &  & j\geq1;
\end{array}
\]%
\[
\pi_{0}=q\pi_{0}+q\pi_{1}
\]
hence%
\[
\pi_{1}=\frac{p}{q}\pi_{0}.
\]
From%
\[
\pi_{1}=p\pi_{0}+q\pi_{2}
\]
we deduce%
\[
q\pi_{2}=\pi_{1}-p\pi_{0}=\left(  1-q\right)  \pi_{1}=p\pi_{1}
\]
hence%
\[
\pi_{2}=\frac{p}{q}\pi_{1}=\frac{p^{2}}{q^{2}}\pi_{0}.
\]
Defining%
\begin{align*}
F(z)  &  =%
{\displaystyle\sum\limits_{j=2}^{\infty}}
\pi_{j}z^{j}=pz%
{\displaystyle\sum\limits_{j=2}^{\infty}}
\pi_{j-1}z^{j-1}+\frac{q}{z}%
{\displaystyle\sum\limits_{j=2}^{\infty}}
\pi_{j+1}z^{j+1}\\
&  =pz\left[  F(z)+\pi_{1}z\right]  +\frac{q}{z}\left[  F(z)-\pi_{2}%
z^{2}\right]  ,
\end{align*}
we have%
\[
\left[  1-pz-\frac{q}{z}\right]  F(z)=p\pi_{1}z^{2}-q\pi_{2}z
\]
hence%
\[
\left[  q-z+pz^{2}\right]  F(z)=\left(  q\pi_{2}-p\pi_{1}z\right)  z
\]
hence%
\[
(1-z)\left(  q-pz\right)  F(z)=\frac{p^{2}\left(  1-z\right)  z}{q}\pi_{0}
\]
hence%
\[
F(z)=\frac{p^{2}z}{q\left(  q-pz\right)  }\pi_{0}
\]
hence%
\[
L=\lim_{z\rightarrow1}F(z)=\frac{p^{2}}{q\left(  q-p\right)  }\pi_{0}.
\]
Because%
\[
\pi_{0}+\pi_{1}+L=1
\]
it follows that%
\[
\left[  1+\frac{p}{q}+\frac{p^{2}}{q\left(  q-p\right)  }\right]  \pi_{0}=1
\]
therefore%
\[
\pi_{0}=\frac{q-p}{q}
\]
and thus%
\[%
\begin{array}
[c]{ccc}%
\pi_{j}=\dfrac{q-p}{q}\left(  \dfrac{p}{q}\right)  ^{j}, &  & j\geq0.
\end{array}
\]
Note that, if $k=\log_{q/p}(n)+h+1$, we have%
\[
\left(  \dfrac{q}{p}\right)  ^{k}=n\left(  \dfrac{q}{p}\right)  ^{h+1}
\]
thus%
\[
\pi_{k}n=\dfrac{q-p}{q}\left(  \dfrac{p}{q}\right)  ^{k}n=\dfrac{q-p}%
{q}\left(  \dfrac{q}{p}\right)  ^{-(h+1)}=\dfrac{p(q-p)}{q^{2}}\left(
\dfrac{q}{p}\right)  ^{-h}.
\]
By \cite{Ald-heu},%
\[
\mathbb{E}(C)=1+p\mathbb{E}(C)+q\left(  \frac{p}{q}\right)  \mathbb{E}(C)
\]
equivalently%
\[
\mathbb{E}(C)=\frac{1}{q-p}
\]
which implies%
\begin{align*}
\mathbb{P}\left\{  M_{n}\leq\log_{q/p}(n)+h\right\}   &  =P\left\{  M_{n}%
<\log_{q/p}(n)+h+1\right\} \\
&  \sim\exp\left[  -\frac{p(q-p)^{2}}{q^{2}}\left(  \frac{q}{p}\right)
^{-h}\right]
\end{align*}
as $n\rightarrow\infty$. \ This formula appears in \cite{Fi2-heu}.
\ A\ similar argument gives an analogous result for random walks with strong
reflection at the origin.

\section{Traffic Light:\ $\ell=1$}

We separate the walk into two subwalks: $i\equiv0\operatorname{mod}2$ and
$i\equiv1\operatorname{mod}2$. \ Let%
\[%
\begin{array}
[c]{ccc}%
U=\left(
\begin{array}
[c]{cccccc}%
q & p & 0 & 0 & 0 & \cdots\\
0 & q & p & 0 & 0 & \cdots\\
0 & 0 & q & p & 0 & \cdots\\
0 & 0 & 0 & q & p & \cdots\\
0 & 0 & 0 & 0 & q & \cdots\\
\vdots & \vdots & \vdots & \vdots & \vdots & \ddots
\end{array}
\right)  , &  & \mathbb{P}\left\{  \left.  j+1\right\vert j\right\}  =p
\end{array}
\]
denote the (infinite) transition matrix from $0$ to $1$, and%
\[%
\begin{array}
[c]{ccc}%
V=\left(
\begin{array}
[c]{cccccc}%
1 & 0 & 0 & 0 & 0 & \cdots\\
q & p & 0 & 0 & 0 & \cdots\\
0 & q & p & 0 & 0 & \cdots\\
0 & 0 & q & p & 0 & \cdots\\
0 & 0 & 0 & q & p & \cdots\\
\vdots & \vdots & \vdots & \vdots & \vdots & \ddots
\end{array}
\right)  , &  & \mathbb{P}\left\{  \left.  j-1\right\vert j\right\}  =q
\end{array}
\]
denote the transition matrix from $1$ to $0$. \ 

\subsection{Subwalk on Even Times}

The subwalk for $i\equiv0\operatorname{mod}2$ has transition matrix%
\[
UV=\left(
\begin{array}
[c]{cccccc}%
(1+p)q & p^{2} & 0 & 0 & 0 & \cdots\\
q^{2} & 2pq & p^{2} & 0 & 0 & \cdots\\
0 & q^{2} & 2pq & p^{2} & 0 & \cdots\\
0 & 0 & q^{2} & 2pq & p^{2} & \cdots\\
0 & 0 & 0 & q^{2} & 2pq & \cdots\\
\vdots & \vdots & \vdots & \vdots & \vdots & \ddots
\end{array}
\right)
\]
thus%
\[%
\begin{array}
[c]{ccc}%
\pi_{j}=p^{2}\pi_{j-1}+2pq\pi_{j}+q^{2}\pi_{j+1}, &  & j\geq1;
\end{array}
\]%
\[
\pi_{0}=(1+p)q\pi_{0}+q^{2}\pi_{1}%
\]
hence%
\[
q^{2}\pi_{1}=\left[  1-(1+p)q\right]  \pi_{0}=\left[  1-\left(  1-p^{2}%
\right)  \right]  \pi_{0}%
\]
hence%
\[
\pi_{1}=\frac{p^{2}}{q^{2}}\pi_{0}.
\]
From%
\[
\pi_{1}=p^{2}\pi_{0}+2pq\pi_{1}+q^{2}\pi_{2}%
\]
we deduce%
\[
q^{2}\pi_{2}=(1-2pq)\pi_{1}-p^{2}\pi_{0}=\left(  1-2pq-q^{2}\right)  \pi
_{1}=p^{2}\pi_{1}%
\]
hence%
\[
\pi_{2}=\frac{p^{2}}{q^{2}}\pi_{1}=\frac{p^{4}}{q^{4}}\pi_{0}.
\]
Defining%
\begin{align*}
F(z) &  =%
{\displaystyle\sum\limits_{j=2}^{\infty}}
\pi_{j}z^{j}\\
&  =p^{2}z%
{\displaystyle\sum\limits_{j=2}^{\infty}}
\pi_{j-1}z^{j-1}+2pq%
{\displaystyle\sum\limits_{j=2}^{\infty}}
\pi_{j}z^{j}+\frac{q^{2}}{z}%
{\displaystyle\sum\limits_{j=2}^{\infty}}
\pi_{j+1}z^{j+1}\\
&  =p^{2}z\left[  F(z)+\pi_{1}z\right]  +2pqF(z)+\frac{q^{2}}{z}\left[
F(z)-\pi_{2}z^{2}\right]
\end{align*}
we have%
\[
\left[  1-p^{2}z-2pq-\frac{q^{2}}{z}\right]  F(z)=p^{2}\pi_{1}z^{2}-q^{2}%
\pi_{2}z
\]
hence%
\[
\left[  q^{2}-(1-2pq)z+p^{2}z^{2}\right]  F(z)=\left(  q^{2}\pi_{2}-p^{2}%
\pi_{1}z\right)  z^{2}%
\]
hence%
\[
(1-z)\left(  q^{2}-p^{2}z\right)  F(z)=\frac{p^{4}\left(  1-z\right)  z^{2}%
}{q^{2}}\pi_{0}%
\]
hence%
\[
F(z)=\frac{p^{4}z^{2}}{q^{2}\left(  q^{2}-p^{2}z\right)  }\pi_{0}%
\]
hence%
\[
L=\lim_{z\rightarrow1}F(z)=\frac{p^{4}}{q^{2}\left(  q^{2}-p^{2}\right)  }%
\pi_{0}=\frac{p^{4}}{q^{2}\left(  q-p\right)  }\pi_{0}.
\]
Because%
\[
\pi_{0}+\pi_{1}+L=1
\]
it follows that%
\[
\left[  1+\frac{p^{2}}{q^{2}}+\frac{p^{4}}{q^{2}\left(  q-p\right)  }\right]
\pi_{0}=1
\]
therefore%
\[
\pi_{0}=\frac{q-p}{q^{2}}%
\]
and thus%
\[%
\begin{array}
[c]{ccc}%
\pi_{j}=\dfrac{q-p}{q^{2}}\left(  \dfrac{p^{2}}{q^{2}}\right)  ^{j}, &  &
j\geq0.
\end{array}
\]
Note that, if $k=\log_{q^{2}/p^{2}}(n)+h+1$, we have%
\[
\left(  \dfrac{q^{2}}{p^{2}}\right)  ^{k}=n\left(  \dfrac{q^{2}}{p^{2}%
}\right)  ^{h+1}%
\]
thus%
\[
\pi_{k}\frac{n}{2}=\dfrac{q-p}{2q^{2}}\left(  \dfrac{p^{2}}{q^{2}}\right)
^{k}n=\dfrac{q-p}{2q^{2}}\left(  \dfrac{q^{2}}{p^{2}}\right)  ^{-(h+1)}%
=\dfrac{p^{2}(q-p)}{2q^{4}}\left(  \dfrac{q^{2}}{p^{2}}\right)  ^{-h}.
\]
By \cite{Ald-heu},%
\[
\mathbb{E}(C)=\frac{1}{q^{2}-p^{2}}=\frac{1}{q-p}%
\]
as before, which implies%
\begin{align*}
\mathbb{P}\left\{  M_{n}\leq\log_{q^{2}/p^{2}}(n)+h\right\}   &  =P\left\{
M_{n}<\log_{q^{2}/p^{2}}(n)+h+1\right\}  \\
&  \sim\exp\left[  -\frac{p^{2}(q-p)^{2}}{2q^{4}}\left(  \frac{q^{2}}{p^{2}%
}\right)  ^{-h}\right]
\end{align*}
as $n\rightarrow\infty$. \ 

\subsection{Subwalk on Odd Times}

The subwalk for $i\equiv1\operatorname{mod}2$ has transition matrix%
\[
VU=\left(
\begin{array}
[c]{cccccc}%
q & p & 0 & 0 & 0 & \cdots\\
q^{2} & 2pq & p^{2} & 0 & 0 & \cdots\\
0 & q^{2} & 2pq & p^{2} & 0 & \cdots\\
0 & 0 & q^{2} & 2pq & p^{2} & \cdots\\
0 & 0 & 0 & q^{2} & 2pq & \cdots\\
\vdots & \vdots & \vdots & \vdots & \vdots & \ddots
\end{array}
\right)
\]
thus%
\[%
\begin{array}
[c]{ccc}%
\pi_{j}=p^{2}\pi_{j-1}+2pq\pi_{j}+q^{2}\pi_{j+1}, &  & j\geq2;
\end{array}
\]%
\[
\pi_{1}=p\pi_{0}+2pq\pi_{1}+q^{2}\pi_{2};
\]%
\[
\pi_{0}=q\pi_{0}+q^{2}\pi_{1}%
\]
hence%
\[
q^{2}\pi_{1}=(1-q)\pi_{0}=p\pi_{0}%
\]
hence%
\[
\pi_{1}=\frac{p}{q^{2}}\pi_{0}.
\]
Also%
\[
q^{2}\pi_{2}=(1-2pq)\pi_{1}-p\pi_{0}=\left(  1-2pq-q^{2}\right)  \pi_{1}%
=p^{2}\pi_{1}%
\]
hence%
\[
\pi_{2}=\frac{p^{2}}{q^{2}}\pi_{1}=\frac{p^{3}}{q^{4}}\pi_{0}.
\]
As earlier, we have%
\[
\left[  q^{2}-(1-2pq)z+p^{2}z^{2}\right]  F(z)=\left(  q^{2}\pi_{2}-p^{2}%
\pi_{1}z\right)  z^{2}%
\]
hence%
\[
(1-z)\left(  q^{2}-p^{2}z\right)  F(z)=\frac{p^{3}\left(  1-z\right)  z^{2}%
}{q^{2}}\pi_{0}%
\]
hence%
\[
F(z)=\frac{p^{3}z^{2}}{q^{2}\left(  q^{2}-p^{2}z\right)  }\pi_{0}%
\]
hence%
\[
L=\lim_{z\rightarrow1}F(z)=\frac{p^{3}}{q^{2}\left(  q^{2}-p^{2}\right)  }%
\pi_{0}=\frac{p^{3}}{q^{2}\left(  q-p\right)  }\pi_{0}.
\]
Because%
\[
\pi_{0}+\pi_{1}+L=1
\]
it follows that%
\[
\left[  1+\frac{p}{q^{2}}+\frac{p^{3}}{q^{2}\left(  q-p\right)  }\right]
\pi_{0}=1
\]
therefore%
\[
\pi_{0}=\frac{q-p}{q}%
\]
and thus%
\[%
\begin{array}
[c]{ccc}%
\pi_{j}=\dfrac{q-p}{pq}\left(  \dfrac{p^{2}}{q^{2}}\right)  ^{j}, &  & j\geq1.
\end{array}
\]
With $k$ as before,
\[
\pi_{k}\frac{n}{2}=\dfrac{q-p}{2pq}\left(  \dfrac{p^{2}}{q^{2}}\right)
^{k}n=\dfrac{q-p}{2pq}\left(  \dfrac{q^{2}}{p^{2}}\right)  ^{-(h+1)}%
=\dfrac{p(q-p)}{2q^{3}}\left(  \dfrac{q^{2}}{p^{2}}\right)  ^{-h}%
\]
and%
\[
\mathbb{P}\left\{  M_{n}\leq\log_{q^{2}/p^{2}}(n)+h\right\}  \sim\exp\left[
-\frac{p(q-p)^{2}}{2q^{3}}\left(  \frac{q^{2}}{p^{2}}\right)  ^{-h}\right]  .
\]
This latter formula is the one we desire and also appears in \cite{Fi2-heu}.
\ Observe that the subwalk exponential coefficient $\varepsilon_{0}$ for
$i\equiv0\operatorname{mod}2$ possesses an extra $p/q$ factor compared to the
coefficient $\varepsilon_{1}$ for $i\equiv1\operatorname{mod}2$,
equivalently,\footnote[1]{Looking ahead, in Section 2, $\varepsilon_{1}$ will
correspond to the subwalk for $i\equiv2\operatorname{mod}4$ and $\varepsilon
_{0}$ will correspond to the subwalk for $i\equiv0\operatorname{mod}4$, that
is, to $V^{2}U^{2}$ and $U^{2}V^{2}$ respectively. \ It can be shown that
$\varepsilon_{0}$ will possess an extra $p^{2}/q^{2}$ factor compared to
$\varepsilon_{1}$. \ In Section 3, $\varepsilon_{1}$ will correspond to the
subwalk for $i\equiv3\operatorname{mod}6$ and $\varepsilon_{0}$ will
correspond to the subwalk for $i\equiv0\operatorname{mod}6$, that is, to
$V^{3}U^{3}$ and $U^{3}V^{3}$ respectively. \ It can be shown that
$\varepsilon_{0}$ will possess an extra $p^{3}/q^{3}$ factor compared to
$\varepsilon_{1}$.}%
\[
\varepsilon_{1}=\frac{1}{2}\cdot\frac{p^{2}}{q^{2}}\cdot\frac{q}{p}\cdot
\frac{(q-p)^{2}}{q^{2}}=\frac{p(q-p)^{2}}{2q^{3}}.
\]
Of course, the two maxima are not independent.

\subsection{Expected Sojourn Time}

In our treatment of the subwalk represented by $UV$, we indicated that
$\mathbb{E}(C)=1/(q-p)$ without comment \cite{Fel-heu}. \ Let us elaborate on
this point. \ Consider a random walk on the integers consisting of incremental
steps satisfying%
\[
\left\{
\begin{array}
[c]{lll}%
-1 &  & \text{with probability }q^{2},\\
0 &  & \text{with probability }2pq,\\
1 &  & \text{with probability }p^{2}.
\end{array}
\right.
\]
\footnote[2]{Defining $\nu_{j}$ is best done as follows. I am at a large
level, say, $J$. \ I place the origin at $J$ and I wish to find the
probability $\nu_{j}$ of returning to $J$ starting from $J-j$, equivalently,
to $0$ now. \ I reverse the walk direction. Now $J-j$ is $j$ and $J+j$ is
$-j$. \ The trend is now towards the positive integers rather than the
negative integers.}For nonzero $j$, let $\nu_{j}$ denote the probability that,
starting from $-j$, the walker eventually hits $0$. \ Let $\nu_{0}$ denote the
probability that, starting from $0$, the walker eventually returns to $0$ (at
some future time). \ We have two values for $\nu_{0}$: when it is used in a
recursion, it is equal to $1$; when it corresponds to a return probability, it
retains the symbol $\nu_{0}$. \ Using%
\[%
\begin{array}
[c]{ccc}%
\nu_{j}=p^{2}\nu_{j-1}+2pq\nu_{j}+q^{2}\nu_{j+1}, &  & j\geq1;
\end{array}
\]%
\[
\nu_{0}=p^{2}\nu_{-1}+2pq+q^{2}\nu_{1}%
\]
define%
\begin{align*}
\tilde{F}(z)  &  =%
{\displaystyle\sum\limits_{j=1}^{\infty}}
\nu_{j}z^{j}\\
&  =p^{2}z%
{\displaystyle\sum\limits_{j=1}^{\infty}}
\nu_{j-1}z^{j-1}+2pq%
{\displaystyle\sum\limits_{j=1}^{\infty}}
\nu_{j}z^{j}+\frac{q^{2}}{z}%
{\displaystyle\sum\limits_{j=1}^{\infty}}
\nu_{j+1}z^{j+1}\\
&  =p^{2}z\left[  \tilde{F}(z)+1\right]  +2pq\tilde{F}(z)+\frac{q^{2}}%
{z}\left[  \tilde{F}(z)-\nu_{1}z\right]
\end{align*}
equivalently%
\[
\left[  1-p^{2}z-2pq-\frac{q^{2}}{z}\right]  \tilde{F}(z)=p^{2}z-q^{2}\nu_{1}%
\]
equivalently%
\[
\left[  q^{2}-(1-2pq)z+p^{2}z^{2}\right]  \tilde{F}(z)=\left(  q^{2}\nu
_{1}-p^{2}z\right)  z
\]
equivalently%
\[
(1-z)\left(  q^{2}-p^{2}z\right)  \tilde{F}(z)=\left(  \nu_{0}-p^{2}\nu
_{-1}-2pq-p^{2}z\right)  z.
\]
Only the first of the zeroes $1$, $q^{2}/p^{2}$ is of interest (the second is
$>1$). \ Substituting $z=1$ into the numerator of $\tilde{F}(z)$ gives an
equation
\[
Eq_{1}:\nu_{0}-p^{2}\nu_{-1}-2pq-p^{2}=0.
\]
Also, using%
\[%
\begin{array}
[c]{ccc}%
\nu_{-j}=p^{2}\nu_{-j-1}+2pq\nu_{-j}+q^{2}\nu_{-j+1}, &  & j\geq1
\end{array}
\]
we deduce that\footnote[3]{The same identity connecting $\nu_{-j}$ and
$\nu_{j}$ will be true in Sections 2 and 3 as well, by the same reasoning.}%
\[
\nu_{-j}=\nu_{j}\left(  \frac{q^{2}}{p^{2}}\right)  ^{j}%
\]
since multiplying both sides of%
\[
\nu_{j}\left(  \frac{q^{2}}{p^{2}}\right)  ^{j}=p^{2}\nu_{j+1}\left(
\frac{q^{2}}{p^{2}}\right)  ^{j+1}+2pq\nu_{j}\left(  \frac{q^{2}}{p^{2}%
}\right)  ^{j}+q^{2}\nu_{j-1}\left(  \frac{q^{2}}{p^{2}}\right)  ^{j-1}%
\]
by $p^{2j}/q^{2j}$ gives an identity. \ Replacing $q^{2}\nu_{1}$ by $p^{2}%
\nu_{-1}$ in our initial expression for $\nu_{0}$ gives another equation
\[
Eq_{2}:\nu_{0}=2p^{2}\nu_{-1}+2pq.
\]
Solving $Eq_{1}$ and $Eq_{2}$ simultaneously yields $\nu_{0}=2p$ and $\nu
_{-1}=1$. \ More generally, $\nu_{-j}=1$ for $j\geq1$. \ Most importantly,
\[
\mathbb{E}(C)=\frac{1}{1-\nu_{0}}=\frac{1}{q-p}%
\]
as was to be shown.

We will similarly study random walks%
\[
\left\{
\begin{array}
[c]{lll}%
-2 &  & \text{with probability }q^{4},\\
-1 &  & \text{with probability }4pq^{3},\\
0 &  & \text{with probability }6p^{2}q^{2},\\
1 &  & \text{with probability }4p^{3}q,\\
2 &  & \text{with probability }p^{4};
\end{array}
\right.
\]%
\[
\left\{
\begin{array}
[c]{lll}%
-3 &  & \text{with probability }q^{6},\\
-2 &  & \text{with probability }6pq^{5},\\
-1 &  & \text{with probability }15p^{2}q^{4}\\
0 &  & \text{with probability }20p^{3}q^{3},\\
1 &  & \text{with probability }15p^{4}q^{2},\\
2 &  & \text{with probability }6p^{5}q,\\
3 &  & \text{with probability }p^{6}%
\end{array}
\right.
\]
in Sections 2.2 and 3.2 respectively. \ The formulas for $\mathbb{E}(C)$,
however, will be somewhat more complicated. \ 

\section{Traffic Light:\ $\ell=2$}

We have four subwalks. \ Let us consider the subwalk represented by
$U^{2}V^{2}$:
\[
\left(
\begin{array}
[c]{ccccccc}%
\left(  1+2p+3p^{2}\right)  q^{2} & 4p^{3}q & p^{4} & 0 & 0 & 0 & \cdots\\
(1+3p)q^{3} & 6p^{2}q^{2} & 4p^{3}q & p^{4} & 0 & 0 & \cdots\\
q^{4} & 4pq^{3} & 6p^{2}q^{2} & 4p^{3}q & p^{4} & 0 & \cdots\\
0 & q^{4} & 4pq^{3} & 6p^{2}q^{2} & 4p^{3}q & p^{4} & \cdots\\
0 & 0 & q^{4} & 4pq^{3} & 6p^{2}q^{2} & 4p^{3}q & \cdots\\
0 & 0 & 0 & q^{4} & 4pq^{3} & 6p^{2}q^{2} & \cdots\\
\vdots & \vdots & \vdots & \vdots & \vdots & \vdots & \ddots
\end{array}
\right)
\]

\subsection{Stationary Distribution}

For $j\geq2$, we have%
\[
\pi_{j}=p^{4}\pi_{j-2}+4p^{3}q\pi_{j-1}+6p^{2}q^{2}\pi_{j}+4pq^{3}\pi
_{j+1}+q^{4}\pi_{j+2};
\]%
\[
\pi_{0}=\left(  1+2p+3p^{2}\right)  q^{2}\pi_{0}+(1+3p)q^{3}\pi_{1}+q^{4}%
\pi_{2}%
\]
hence%
\[
\pi_{2}=\frac{\left(  1-q^{2}-2pq^{2}-3p^{2}q^{2}\right)  \pi_{0}%
-(1+3p)q^{3}\pi_{1}}{q^{4}};
\]%
\[
\pi_{1}=4p^{3}q\pi_{0}+6p^{2}q^{2}\pi_{1}+4pq^{3}\pi_{2}+q^{4}\pi_{3}%
\]
hence%
\[
\pi_{3}=\frac{-4p^{3}q\pi_{0}+\left(  1-6p^{2}q^{2}\right)  \pi_{1}-4pq^{3}%
\pi_{2}}{q^{4}}.
\]
Defining%
\begin{align*}
F(z)  &  =%
{\displaystyle\sum\limits_{j=2}^{\infty}}
\pi_{j}z^{j}\\
&  =p^{4}z^{2}%
{\displaystyle\sum\limits_{j=2}^{\infty}}
\pi_{j-2}z^{j-2}+4p^{3}qz%
{\displaystyle\sum\limits_{j=2}^{\infty}}
\pi_{j-1}z^{j-1}+6p^{2}q^{2}%
{\displaystyle\sum\limits_{j=2}^{\infty}}
\pi_{j}z^{j}\\
&  +\frac{4pq^{3}}{z}%
{\displaystyle\sum\limits_{j=2}^{\infty}}
\pi_{j+1}z^{j+1}+\frac{q^{4}}{z^{2}}%
{\displaystyle\sum\limits_{j=2}^{\infty}}
\pi_{j+2}z^{j+2}\\
&  =p^{4}z^{2}\left[  F(z)+\pi_{0}+\pi_{1}z\right]  +4p^{3}qz\left[
F(z)+\pi_{1}z\right]  +6p^{2}q^{2}F(z)\\
&  +\frac{4pq^{3}}{z}\left[  F(z)-\pi_{2}z^{2}\right]  +\frac{q^{4}}{z^{2}%
}\left[  F(z)-\pi_{2}z^{2}-\pi_{3}z^{3}\right]
\end{align*}
we deduce%
\begin{align*}
&  \left[  1-p^{4}z^{2}-4p^{3}qz-6p^{2}q^{2}-\frac{4pq^{3}}{z}-\frac{q^{4}%
}{z^{2}}\right]  F(z)\\
&  =p^{4}z^{2}\left(  \pi_{0}+\pi_{1}z\right)  +4p^{3}qz\left(  \pi
_{1}z\right)  -\frac{4pq^{3}}{z}\left(  \pi_{2}z^{2}\right)  -\frac{q^{4}%
}{z^{2}}\left(  \pi_{2}z^{2}+\pi_{3}z^{3}\right)
\end{align*}
hence%
\begin{align*}
&  \left[  q^{4}+4pq^{3}z-\left(  1-6p^{2}q^{2}\right)  z^{2}+4p^{3}%
qz^{3}+p^{4}z^{4}\right]  F(z)\\
&  =-p^{4}z^{4}\left(  \pi_{0}+\pi_{1}z\right)  -4p^{3}qz^{3}\left(  \pi
_{1}z\right)  +4pq^{3}z\left(  \pi_{2}z^{2}\right)  +q^{4}\left(  \pi_{2}%
z^{2}+\pi_{3}z^{3}\right)  .
\end{align*}
Replacing $\pi_{2}$ and $\pi_{3}$ by expressions in $\pi_{0}$ and $\pi_{1}$,
then cancelling the common factor $1-z$ between numerator and denominator,
yields
\[
F(z)=\frac{\left\{  p^{3}(4-3p+pz)\pi_{0}+\left[  -1+6p^{2}-8p^{3}%
+3p^{4}+(4-3p)p^{3}z+p^{4}z^{2}\right]  \pi_{1}\right\}  z^{2}}{(q^{2}%
-p^{2}z)\left[  q^{2}+(1+2pq)z+p^{2}z^{2}\right]  }%
\]
hence%
\[
L=\lim_{z\rightarrow1}F(z)=\frac{2p^{3}(1+q)\pi_{0}-\left[  2(q-p)-q^{4}%
\right]  \pi_{1}}{2(q-p)}.
\]
We observe three zeroes in the denominator $D(z)$ of $F(z)$. \ The first zero,
of smallest modulus $<1$, is negative and given by%
\[
z_{1}=\frac{-1-2pq+\theta}{2p^{2}}%
\]
where $\theta=\sqrt{1+4pq}$. \ The second zero, of intermediate modulus, is
positive and given by%
\[
z_{2}=\frac{q^{2}}{p^{2}}>1.
\]
The third zero, of largest modulus $>1$, is negative and given by%
\[
z_{3}=\frac{-1-2pq-\theta}{2p^{2}}.
\]
Finding the unknowns $\pi_{0}$ and $\pi_{1}$ is achieved by solving two
simultaneous equations:%
\[
Eq_{1}:subst\left(  z=z_{1},N\right)  =0
\]
(substituting $z_{1}$ for $z$ in the numerator $N(z)$ for $F(z)$ and setting
this equal to zero)
\[
Eq_{2}:\pi_{0}+\pi_{1}+L=1
\]
which yields%
\[
\pi_{0}=\frac{(q-p)\left(  3-2p-\theta\right)  }{2q^{4}},
\]%
\[
\pi_{1}=\frac{(q-p)\left[  -1-p-2pq+(1+p)\theta\right]  }{q^{5}}.
\]
Thus we have a complete description of the stationary distribution. \ An exact
expression for $\pi_{j}$ is infeasible; therefore asymptotics as
$j\rightarrow\infty$ are necessary. \ The second zero $z_{2}$ leads, by
classical singularity analysis, to \cite{FS-heu}%
\[
A(p)=-\frac{N(z_{2})}{z_{2}D^{\prime}(z_{2})}=\frac{(q-p)\left[
1+(q-p)\theta\right]  }{4q^{4}},
\]%
\[
\pi_{j}\sim A(p)\left(  \dfrac{p^{2}}{q^{2}}\right)  ^{j}.
\]
This is the expression that we shall use in the clumping heuristic.

\subsection{Clump Rate}

Using%
\[%
\begin{array}
[c]{ccc}%
\nu_{j}=p^{4}\nu_{j-2}+4p^{3}q\nu_{j-1}+6p^{2}q^{2}\nu_{j}+4pq^{3}\nu
_{j+1}+q^{4}\nu_{j+2}, &  & j\geq1;
\end{array}
\]%
\[
\nu_{0}=p^{4}\nu_{-2}+4p^{3}q\nu_{-1}+6p^{2}q^{2}+4pq^{3}\nu_{1}+q^{4}\nu_{2}
\]
define%
\begin{align*}
\tilde{F}(z)  &  =%
{\displaystyle\sum\limits_{j=1}^{\infty}}
\nu_{j}z^{j}\\
&  =p^{4}z^{2}%
{\displaystyle\sum\limits_{j=1}^{\infty}}
\nu_{j-2}z^{j-2}+4p^{3}qz%
{\displaystyle\sum\limits_{j=1}^{\infty}}
\nu_{j-1}z^{j-1}+6p^{2}q^{2}%
{\displaystyle\sum\limits_{j=1}^{\infty}}
\nu_{j}z^{j}\\
&  +\frac{4pq^{3}}{z}%
{\displaystyle\sum\limits_{j=1}^{\infty}}
\nu_{j+1}z^{j+1}+\frac{q^{4}}{z^{2}}%
{\displaystyle\sum\limits_{j=1}^{\infty}}
\nu_{j+2}z^{j+2}\\
&  =p^{4}z^{2}\left[  \tilde{F}(z)+\nu_{-1}z^{-1}+1\right]  +4p^{3}qz\left[
\tilde{F}(z)+1\right]  +6p^{2}q^{2}\tilde{F}(z)\\
&  +\frac{4pq^{3}}{z}\left[  \tilde{F}(z)-\nu_{1}z\right]  +\frac{q^{4}}%
{z^{2}}\left[  \tilde{F}(z)-\nu_{1}z-\nu_{2}z^{2}\right]
\end{align*}
equivalently%
\begin{align*}
&  \left[  1-p^{4}z^{2}-4p^{3}qz-6p^{2}q^{2}-\frac{4pq^{3}}{z}-\frac{q^{4}%
}{z^{2}}\right]  \tilde{F}(z)\\
&  =p^{4}z^{2}\left(  \nu_{-1}z^{-1}+1\right)  +4p^{3}qz-\frac{4pq^{3}}%
{z}\left(  \nu_{1}z\right)  -\frac{q^{4}}{z^{2}}\left(  \nu_{1}z+\nu_{2}%
z^{2}\right)
\end{align*}
equivalently%
\begin{align*}
&  \left[  q^{4}+4pq^{3}z-\left(  1-6p^{2}q^{2}\right)  z^{2}+4p^{3}%
qz^{3}+p^{4}z^{4}\right]  \tilde{F}(z)\\
&  =-p^{4}z^{4}\left(  \nu_{-1}z^{-1}+1\right)  -4p^{3}qz^{3}+4pq^{3}z\left(
\nu_{1}z\right)  +q^{4}\left(  \nu_{1}z+\nu_{2}z^{2}\right)
\end{align*}
equivalently%
\begin{align*}
&  (1-z)(q^{2}-p^{2}z)\left[  q^{2}+(1+2pq)z+p^{2}z^{2}\right]  \tilde{F}(z)\\
&  =-p^{4}z^{3}\nu_{-1}-p^{4}z^{4}-4p^{3}qz^{3}+4pq^{3}z^{2}\nu_{1}+q^{4}%
z\nu_{1}\\
&  +z^{2}\left(  \nu_{0}-p^{4}\nu_{-2}-4p^{3}q\nu_{-1}-6p^{2}q^{2}-4pq^{3}%
\nu_{1}\right) \\
&  =z^{2}\nu_{0}+q^{4}z\nu_{1}-p^{4}z^{3}\nu_{-1}-4p^{3}qz^{2}\nu_{-1}%
-p^{4}z^{2}\nu_{-2}-6p^{2}q^{2}z^{2}-4p^{3}qz^{3}-p^{4}z^{4}.
\end{align*}
Only the first two of the four zeroes $z_{1}$, $1$, $z_{2}$, $z_{3}$ are of
interest. \ Let $\tilde{N}(z)$ denote the numerator for $\tilde{F}(z)$. \ We
have%
\[
\tilde{E}q_{1}:subst\left(  z=z_{1},\tilde{N}\right)  =0,
\]%
\[
\tilde{E}q_{2}:subst\left(  z=1,\tilde{N}\right)  =0.
\]
Replacing $q^{4}\nu_{2}$ by $p^{4}\nu_{-2}$ in our initial expression for
$\nu_{0}$ gives%
\[
\tilde{E}q_{3}:\nu_{0}=2p^{4}\nu_{-2}+4p^{3}q\nu_{-1}+6p^{2}q^{2}+4pq^{3}%
\nu_{1}.
\]
Also, replacing $q^{2}\nu_{1}$ by $p^{2}\nu_{-1}$ throughout $\tilde{E}q_{1}$,
$\tilde{E}q_{2}$ and $\tilde{E}q_{3}$ reduces the number of variables to
three. \ The simultaneous solution is
\[
\nu_{0}=\frac{-1+2p+8p^{2}-8p^{3}+(q-p)^{2}\theta}{4pq},
\]%
\[
\nu_{-1}=\frac{1-8p^{2}+16p^{3}-8p^{4}-(q-p)\theta}{8p^{3}q},
\]%
\[
\nu_{-2}=\frac{-1-2p+12p^{2}-24p^{4}+24p^{5}-8p^{6}+(q-p)(1+2p-4p^{2})\theta
}{8p^{5}q}
\]
yielding%
\[
\nu_{1}=\frac{1-8p^{2}+16p^{3}-8p^{4}-(q-p)\theta}{8pq^{3}}
\]
in particular. \ 

\footnote[4]{Employing the original coordinate axis may aid understanding.
\ Readers might be tempted to use the level $j$ as the absorbing set $S$. But
the maximum could be above $j$ without ever touching level $j$ because of the
transition $p^{4}$. So we must use as $S$ the levels $j$ and $j+1$: no maximum
can be above $j+1$ without touching at least one of the levels $j$ or $j+1$.
In the revised notation, this leads to the absorbing set $\Omega=\{0,-1\}$%
.}Readers might be tempted to use $\{0\}$ as the absorbing set $\Omega$,
imitating what we did in Section 1.3. \ But, starting from a negative integer,
the walker could stray into the positive integers without ever touching $0$
due to the transition $p^{4}$. \ We must take $\Omega=\{0,-1\}$: no walk can
venture above $0$ without touching at least one of $0$ or $-1$.\ 

An idea of Aldous \cite{Ald-heu} now comes crucially into play. \ The rate
$\lambda$ of clumps of visits to $\Omega$ is equal to $\lambda_{0}%
+\lambda_{-1}$ where parameters $\lambda_{0}$ and $\lambda_{-1}$ are solutions
of the system%
\[
\lambda_{0}+\lambda_{-1}\nu_{-1}=\left(  1-\nu_{0}\right)  \pi_{j},
\]%
\[
\lambda_{0}\nu_{1}+\lambda_{-1}=\left(  1-\nu_{0}\right)  \pi_{j+1}\sim
\frac{p^{2}}{q^{2}}\left(  1-\nu_{0}\right)  \pi_{j}.
\]
In words, for nonzero $j$, the ratio $\nu_{j}/\left(  1-\nu_{0}\right)  $ is
the expected sojourn time in $\{0\}$, given that the walk started at $-j$.
\ The total clump rate is consequently%
\[
\lambda\sim\frac{(q-p)\left[  1+(q-p)\theta\right]  }{2q^{2}}\pi_{j}%
\]
and thus the exponential coefficient is%
\[
\varepsilon_{1}=\frac{1}{4}\cdot\frac{p^{2}}{q^{2}}\cdot\frac{q^{2}}{p^{2}%
}\cdot\frac{\lambda}{\pi_{j}}\cdot A(p)\sim\frac{(q-p)^{2}\left[
1+(q-p)\theta\right]  ^{2}}{32q^{6}}=\frac{\chi_{2}(p)}{4}%
\]
where%
\[
\chi_{2}(p)=\frac{(q-p)^{2}}{4q^{6}}\left[  \left(  1-8p^{2}+16p^{3}%
-8p^{4}\right)  +(q-p)\theta\right]
\]
was conjectured in \cite{Fi1-heu}. \ What was missed previously, however, is
the expression for $\varepsilon_{0}$ as a perfect square. \ This fact is a
corollary of the hidden relation
\[
\frac{\lambda}{\pi_{j}}=2q^{2}A(p)
\]
(\textquotedblleft hidden\textquotedblright\ in the sense that our
experimental methods in \cite{Fi1-heu} overlooked this intriguing formula).

\section{Traffic Light:\ $\ell=3$}

We have six subwalks. \ Let us consider the subwalk represented by $U^{3}%
V^{3}$:
\[
\left(  {\small
\begin{array}
[c]{ccccccccc}%
\left(  1+3p+6p^{2}+10p^{3}\right)  q^{3} & 15p^{4}q^{2} & 6p^{5}q & p^{6} &
0 & 0 & 0 & 0 & \cdots\\
\left(  1+4p+10p^{2}\right)  q^{4} & 20p^{3}q^{3} & 15p^{4}q^{2} & 6p^{5}q &
p^{6} & 0 & 0 & 0 & \cdots\\
(1+5p)q^{5} & 15p^{2}q^{4} & 20p^{3}q^{3} & 15p^{4}q^{2} & 6p^{5}q & p^{6} &
0 & 0 & \cdots\\
q^{6} & 6pq^{5} & 15p^{2}q^{4} & 20p^{3}q^{3} & 15p^{4}q^{2} & 6p^{5}q & p^{6}
& 0 & \cdots\\
0 & q^{6} & 6pq^{5} & 15p^{2}q^{4} & 20p^{3}q^{3} & 15p^{4}q^{2} & 6p^{5}q &
p^{6} & \cdots\\
0 & 0 & q^{6} & 6pq^{5} & 15p^{2}q^{4} & 20p^{3}q^{3} & 15p^{4}q^{2} &
6p^{5}q & \cdots\\
0 & 0 & 0 & q^{6} & 6pq^{5} & 15p^{2}q^{4} & 20p^{3}q^{3} & 15p^{4}q^{2} &
\cdots\\
0 & 0 & 0 & 0 & q^{6} & 6pq^{5} & 15p^{2}q^{4} & 20p^{3}q^{3} & \cdots\\
\vdots & \vdots & \vdots & \vdots & \vdots & \vdots & \vdots & \vdots & \ddots
\end{array}
}\right)
\]

\subsection{Stationary Distribution}

For $j\geq3$, we have%
\[
\pi_{j}=p^{6}\pi_{j-3}+6p^{5}q\pi_{j-2}+15p^{4}q^{2}\pi_{j-1}+20p^{3}q^{3}%
\pi_{j}+15p^{2}q^{4}\pi_{j+1}+6pq^{5}\pi_{j+2}+q^{6}\pi_{j+3};
\]%
\[
\pi_{0}=\left(  1+3p+6p^{2}+10p^{3}\right)  q^{3}\pi_{0}+\left(
1+4p+10p^{2}\right)  q^{4}\pi_{1}+(1+5p)q^{5}\pi_{2}+q^{6}\pi_{3}
\]
hence%
\[
\pi_{3}=\frac{\left(  1-q^{3}-3pq^{3}-6p^{2}q^{3}-10p^{3}q^{3}\right)  \pi
_{0}-(1+4p+10p^{2})q^{4}\pi_{1}-(1+5p)q^{5}\pi_{2}}{q^{6}};
\]%
\[
\pi_{1}=15p^{4}q^{2}\pi_{0}+20p^{3}q^{3}\pi_{1}+15p^{2}q^{4}\pi_{2}+6pq^{5}%
\pi_{3}+q^{6}\pi_{4}
\]
hence%
\[
\pi_{4}=\frac{-15p^{4}q^{2}\pi_{0}+\left(  1-20p^{3}q^{3}\right)  \pi
_{1}-15p^{2}q^{4}\pi_{2}-6pq^{5}\pi_{3}}{q^{6}};
\]%
\[
\pi_{2}=6p^{5}q\pi_{0}+15p^{4}q^{2}\pi_{1}+20p^{3}q^{3}\pi_{2}+15p^{2}q^{4}%
\pi_{3}+6pq^{5}\pi_{4}+q^{6}\pi_{5}
\]
hence%
\[
\pi_{5}=\frac{-6p^{5}q\pi_{0}-15p^{4}q^{2}\pi_{1}+\left(  1-20p^{3}%
q^{3}\right)  \pi_{2}-15p^{2}q^{4}\pi_{3}-6pq^{5}\pi_{4}}{q^{6}}.
\]
Defining%
\begin{align*}
F(z)  &  =%
{\displaystyle\sum\limits_{j=3}^{\infty}}
\pi_{j}z^{j}\\
&  =p^{6}z^{3}%
{\displaystyle\sum\limits_{j=3}^{\infty}}
\pi_{j-3}z^{j-3}+6p^{5}qz^{2}%
{\displaystyle\sum\limits_{j=3}^{\infty}}
\pi_{j-2}z^{j-2}\\
&  +15p^{4}q^{2}z%
{\displaystyle\sum\limits_{j=3}^{\infty}}
\pi_{j-1}z^{j-1}+20p^{3}q^{3}%
{\displaystyle\sum\limits_{j=3}^{\infty}}
\pi_{j}z^{j}+\frac{15p^{2}q^{4}}{z}%
{\displaystyle\sum\limits_{j=3}^{\infty}}
\pi_{j+1}z^{j+1}\\
&  +\frac{6pq^{5}}{z^{2}}%
{\displaystyle\sum\limits_{j=3}^{\infty}}
\pi_{j+2}z^{j+2}+\frac{q^{6}}{z^{3}}%
{\displaystyle\sum\limits_{j=3}^{\infty}}
\pi_{j+3}z^{j+3}\\
&  =p^{6}z^{3}\left[  F(z)+\pi_{0}+\pi_{1}z+\pi_{2}z^{2}\right]  +6p^{5}%
qz^{2}\left[  F(z)+\pi_{1}z+\pi_{2}z^{2}\right] \\
&  +15p^{4}q^{2}z\left[  F(z)+\pi_{2}z^{2}\right]  +20p^{3}q^{3}%
F(z)+\frac{15p^{2}q^{4}}{z}\left[  F(z)-\pi_{3}z^{3}\right] \\
&  +\frac{6pq^{5}}{z^{2}}\left[  F(z)-\pi_{3}z^{3}-\pi_{4}z^{4}\right]
+\frac{q^{6}}{z^{3}}\left[  F(z)-\pi_{3}z^{3}-\pi_{4}z^{4}-\pi_{5}%
z^{5}\right]
\end{align*}
we deduce%
\begin{align*}
&  \left[  1-p^{6}z^{3}-6p^{5}qz^{2}-15p^{4}q^{2}z-20p^{3}q^{3}-\frac
{15p^{2}q^{4}}{z}-\frac{6pq^{5}}{z^{2}}-\frac{q^{6}}{z^{3}}\right]  F(z)\\
&  =p^{6}z^{3}\left(  \pi_{0}+\pi_{1}z+\pi_{2}z^{2}\right)  +6p^{5}%
qz^{2}\left(  \pi_{1}z+\pi_{2}z^{2}\right)  +15p^{4}q^{2}z\left(  \pi_{2}%
z^{2}\right) \\
&  -\frac{15p^{2}q^{4}}{z}\left(  \pi_{3}z^{3}\right)  -\frac{6pq^{5}}{z^{2}%
}\left(  \pi_{3}z^{3}+\pi_{4}z^{4}\right)  -\frac{q^{6}}{z^{3}}\left(  \pi
_{3}z^{3}+\pi_{4}z^{4}+\pi_{5}z^{5}\right)
\end{align*}
hence%
\begin{align*}
&  \left[  q^{6}+6pq^{5}z+15p^{2}q^{4}z^{2}-\left(  1-20p^{3}q^{3}\right)
z^{3}+15p^{4}q^{2}z^{4}+6p^{5}qz^{5}+p^{6}z^{6}\right]  F(z)\\
&  =-p^{6}z^{6}\left(  \pi_{0}+\pi_{1}z+\pi_{2}z^{2}\right)  -6p^{5}%
qz^{5}\left(  \pi_{1}z+\pi_{2}z^{2}\right)  -15p^{4}q^{2}z^{4}\left(  \pi
_{2}z^{2}\right) \\
&  +15p^{2}q^{4}z^{2}\left(  \pi_{3}z^{3}\right)  +6pq^{5}z\left(  \pi
_{3}z^{3}+\pi_{4}z^{4}\right)  +q^{6}\left(  \pi_{3}z^{3}+\pi_{4}z^{4}+\pi
_{5}z^{5}\right)  .
\end{align*}
Replacing $\pi_{3}$, $\pi_{4}$ and $\pi_{5}$ by expressions in $\pi_{0}$,
$\pi_{1}$ and $\pi_{2}$, then cancelling the common factor $1-z$ between
numerator and denominator, yields $F(z)$ to be%
\[
\frac{\left\{  p^{4}\left[  b+apz+p^{2}z^{2}\right]  \pi_{0}+\left[
c+bp^{4}z+ap^{5}z^{2}+p^{6}z^{3}\right]  \pi_{1}+\left[  d+cz+bp^{4}%
z^{2}+ap^{5}z^{3}+p^{6}z^{4}\right]  \pi_{2}\right\}  z^{3}}{(q^{2}%
-p^{2}z)\left[  q^{4}+q^{2}(1+4pq)z+(1+2pq+6p^{2}q^{2})z^{2}+p^{2}%
(1+4pq)z^{3}+p^{4}z^{4}\right]  }
\]
where%
\[%
\begin{array}
[c]{ccc}%
a=6-5p, &  & b=15-24p+10p^{2},
\end{array}
\]%
\[%
\begin{array}
[c]{ccc}%
c=-\left(  1-20p^{3}+45p^{4}-36p^{5}+10p^{6}\right)  , &  & d=-\left(
1-15p^{2}+40p^{3}-45p^{4}+24p^{5}-5p^{6}\right)
\end{array}
\]
hence%
\[
L=\lim_{z\rightarrow1}F(z)=\frac{3p^{4}\left(  1+2q+2q^{2}\right)  \pi
_{0}-\left[  3(q-p)-2(1+2p)q^{5}\right]  \pi_{1}-\left[  3(q-p)-q^{6}\right]
\pi_{2}}{3(q-p)}.
\]
We observe five zeroes in the denominator $D(z)$ of $F(z)$. \ Two (complex
conjugate) zeroes have modulus $<1$:%
\[
z_{1}=\frac{-1-i\sqrt{3}-4pq+\sqrt{-2+2i\sqrt{3}+8\left(  1+i\sqrt{3}\right)
pq}}{4p^{2}},
\]%
\[
z_{2}=\frac{-1+i\sqrt{3}-4pq+\sqrt{-2-2i\sqrt{3}+8\left(  1-i\sqrt{3}\right)
pq}}{4p^{2}};
\]
two zeroes have modulus $>1$:%
\[
z_{4}=\frac{-1-i\sqrt{3}-4pq-\sqrt{-2+2i\sqrt{3}+8\left(  1+i\sqrt{3}\right)
pq}}{4p^{2}},
\]%
\[
z_{5}=\frac{-1+i\sqrt{3}-4pq-\sqrt{-2-2i\sqrt{3}+8\left(  1-i\sqrt{3}\right)
pq}}{4p^{2}};
\]
and the remaining (real) zero is%
\[
z_{3}=\frac{q^{2}}{p^{2}}>1.
\]
Finding the unknowns $\pi_{0}$, $\pi_{1}$ and $\pi_{2}$ is achieved by solving
three simultaneous equations (two involving the numerator $N(z)$ of $F(z)$):%
\[
Eq_{1}:subst\left(  z=z_{1},N\right)  =0
\]%
\[
Eq_{2}:subst\left(  z=z_{2},N\right)  =0
\]%
\[
Eq_{3}:\pi_{0}+\pi_{1}+\pi_{2}+L=1
\]
which yields
\[
\pi_{0}=\tfrac{(q-p)\left[  7-10p+4p^{2}+\theta-\sqrt{2}\sqrt{1+28p-60p^{2}%
+40p^{3}-8p^{4}+\left(  7-10p+4p^{2}\right)  \theta}\right]  }{4q^{6}},
\]%
\begin{align*}
\pi_{1}  &  =\tfrac{(q-p)\left[  -3\left(  1+7p-10p^{2}+4p^{3}\right)
-3(1+p)\theta+\sqrt{6}\right.  }{{}}\\
&  \tfrac{\left.  \cdot\sqrt{-1+30p+71p^{2}-84p^{3}-100p^{4}+120p^{5}%
-24p^{6}+\left(  7+16p-5p^{2}-18p^{3}+12p^{4}\right)  \theta}\right]  }%
{4q^{7}},
\end{align*}%
\begin{align*}
\pi_{2}  &  =\tfrac{(q-p)\left[  3\left(  -1+6p+14p^{2}-20p^{3}+8p^{4}\right)
+3\left(  1+4p+2p^{2}\right)  \theta-\sqrt{6}\right.  }{{}}\\
&  \tfrac{\left.  \cdot\sqrt{-1-16p+64p^{2}+656p^{3}+52p^{4}-1072p^{5}%
+80p^{6}+480p^{7}-96p^{8}+\left(  1+14p+120p^{2}+80p^{3}-92p^{4}%
-24p^{5}+48p^{6}\right)  \theta}\right]  }{4q^{8}}%
\end{align*}
where $\theta=\sqrt{1+4pq+16p^{2}q^{2}}$. \ Thus we have a complete
description of the stationary distribution. \ An exact expression for $\pi
_{j}$ is again infeasible. \ The third zero $z_{3}$ leads to
\[
A(p)=-\frac{N(z_{3})}{z_{3}D^{\prime}(z_{3})}=\frac{(q-p)u+(q-p)^{3}%
\theta+\sqrt{2}(q-p)^{2}\sqrt{v+u\theta}}{12q^{6}},
\]%
\[
\pi_{j}\sim A(p)\left(  \dfrac{p^{2}}{q^{2}}\right)  ^{j}
\]
where%
\[%
\begin{array}
[c]{ccc}%
u=1-2p+6p^{2}-8p^{3}+4p^{4}, &  & v=1+6p^{2}-28p^{3}+54p^{4}-48p^{5}+16p^{6}.
\end{array}
\]
This is the expression that we shall use in the clumping heuristic.

\subsection{Clump Rate}

Using%
\[%
\begin{array}
[c]{ccc}%
\nu_{j}=p^{6}\nu_{j-3}+6p^{5}q\nu_{j-2}+15p^{4}q^{2}\nu_{j-1}+20p^{3}q^{3}%
\nu_{j}+15p^{2}q^{4}\nu_{j+1}+6pq^{5}\nu_{j+2}+q^{6}\nu_{j+3}, &  & j\geq1;
\end{array}
\]%
\[
\nu_{0}=p^{6}\nu_{-3}+6p^{5}q\nu_{-2}+15p^{4}q^{2}\nu_{-1}+20p^{3}%
q^{3}+15p^{2}q^{4}\nu_{1}+6pq^{5}\nu_{2}+q^{6}\nu_{3}
\]
define%
\begin{align*}
\tilde{F}(z)  &  =%
{\displaystyle\sum\limits_{j=1}^{\infty}}
\nu_{j}z^{j}\\
&  =p^{6}z^{3}%
{\displaystyle\sum\limits_{j=1}^{\infty}}
\nu_{j-3}z^{j-3}+6p^{5}qz^{2}%
{\displaystyle\sum\limits_{j=1}^{\infty}}
\nu_{j-2}z^{j-2}\\
&  +15p^{4}q^{2}z%
{\displaystyle\sum\limits_{j=1}^{\infty}}
\nu_{j-1}z^{j-1}+20p^{3}q^{3}%
{\displaystyle\sum\limits_{j=1}^{\infty}}
\nu_{j}z^{j}+\frac{15p^{2}q^{4}}{z}%
{\displaystyle\sum\limits_{j=1}^{\infty}}
\nu_{j+1}z^{j+1}\\
&  +\frac{6pq^{5}}{z^{2}}%
{\displaystyle\sum\limits_{j=1}^{\infty}}
\nu_{j+2}z^{j+2}+\frac{q^{6}}{z^{3}}%
{\displaystyle\sum\limits_{j=1}^{\infty}}
\nu_{j+3}z^{j+3}\\
&  =p^{6}z^{3}\left[  \tilde{F}(z)+\nu_{-2}z^{-2}+\nu_{-1}z^{-1}+1\right]
+6p^{5}qz^{2}\left[  \tilde{F}(z)+\nu_{-1}z^{-1}+1\right] \\
&  +15p^{4}q^{2}z\left[  \tilde{F}(z)+1\right]  +20p^{3}q^{3}\tilde
{F}(z)+\frac{15p^{2}q^{4}}{z}\left[  \tilde{F}(z)-\nu_{1}z\right] \\
&  +\frac{6pq^{5}}{z^{2}}\left[  \tilde{F}(z)-\nu_{1}z-\nu_{2}z^{2}\right]
+\frac{q^{6}}{z^{3}}\left[  \tilde{F}(z)-\nu_{1}z-\nu_{2}z^{2}-\nu_{3}%
z^{3}\right]
\end{align*}
equivalently%
\begin{align*}
&  \left[  1-p^{6}z^{3}-6p^{5}qz^{2}-15p^{4}q^{2}z-20p^{3}q^{3}-\frac
{15p^{2}q^{4}}{z}-\frac{6pq^{5}}{z^{2}}-\frac{q^{6}}{z^{3}}\right]  \tilde
{F}(z)\\
&  =p^{6}z^{3}\left(  \nu_{-2}z^{-2}+\nu_{-1}z^{-1}+1\right)  +6p^{5}%
qz^{2}\left(  \nu_{-1}z^{-1}+1\right)  +15p^{4}q^{2}z\\
&  -\frac{15p^{2}q^{4}}{z}\left(  \nu_{1}z\right)  -\frac{6pq^{5}}{z^{2}%
}\left(  \nu_{1}z+\nu_{2}z^{2}\right)  -\frac{q^{6}}{z^{3}}\left(  \nu
_{1}z+\nu_{2}z^{2}+\nu_{3}z^{3}\right)
\end{align*}
equivalently%
\begin{align*}
&  \left[  q^{6}+6pq^{5}z+15p^{2}q^{4}z^{2}-\left(  1-20p^{3}q^{3}\right)
z^{3}+15p^{4}q^{2}z^{4}+6p^{5}qz^{5}+p^{6}z^{6}\right]  \tilde{F}(z)\\
&  =-p^{6}z^{6}\left(  \nu_{-2}z^{-2}+\nu_{-1}z^{-1}+1\right)  -6p^{5}%
qz^{5}\left(  \nu_{-1}z^{-1}+1\right)  -15p^{4}q^{2}z^{4}\\
&  +15p^{2}q^{4}z^{2}\left(  \nu_{1}z\right)  +6pq^{5}z\left(  \nu_{1}%
z+\nu_{2}z^{2}\right)  +q^{6}\left(  \nu_{1}z+\nu_{2}z^{2}+\nu_{3}%
z^{3}\right)
\end{align*}
equivalently
\begin{align*}
&  (1-z)(q^{2}-p^{2}z)\left[  q^{4}+q^{2}(1+4pq)z+(1+2pq+6p^{2}q^{2}%
)z^{2}+p^{2}(1+4pq)z^{3}+p^{4}z^{4}\right]  \tilde{F}(z)\\
&  =-p^{6}z^{4}\nu_{-2}-p^{6}z^{5}\nu_{-1}-p^{6}z^{6}-6p^{5}qz^{4}\nu
_{-1}-6p^{5}qz^{5}-15p^{4}q^{2}z^{4}\\
&  +15p^{2}q^{4}z^{3}\nu_{1}+6pq^{5}z^{2}\nu_{1}+6pq^{5}z^{3}\nu_{2}+q^{6}%
z\nu_{1}+q^{6}z^{2}\nu_{2}\\
&  +z^{3}\left(  \nu_{0}-p^{6}\nu_{-3}-6p^{5}q\nu_{-2}-15p^{4}q^{2}\nu
_{-1}-20p^{3}q^{3}-15p^{2}q^{4}\nu_{1}-6pq^{5}\nu_{2}\right) \\
&  =z^{3}\nu_{0}+q^{6}z\nu_{1}+6pq^{5}z^{2}\nu_{1}+q^{6}z^{2}\nu_{2}%
-15p^{4}q^{2}z^{3}\nu_{-1}-6p^{5}qz^{4}\nu_{-1}-p^{6}z^{5}\nu_{-1}\\
&  -6p^{5}qz^{3}\nu_{-2}-p^{6}z^{4}\nu_{-2}-p^{6}z^{3}\nu_{-3}-20p^{3}%
q^{3}z^{3}-15p^{4}q^{2}z^{4}-6p^{5}qz^{5}-p^{6}z^{6}.
\end{align*}
Only the first three of the six zeroes $z_{1}$, $z_{2}$, $1$, $z_{3}$, $z_{4}
$, $z_{5}$ are of interest. \ Let $\tilde{N}(z)$ denote the numerator for
$\tilde{F}(z)$. \ We have%
\[
\tilde{E}q_{1}:subst\left(  z=z_{1},\tilde{N}\right)  =0,
\]%
\[
\tilde{E}q_{2}:subst\left(  z=z_{2},\tilde{N}\right)  =0,
\]%
\[
\tilde{E}q_{3}:subst\left(  z=1,\tilde{N}\right)  =0.
\]
Replacing $q^{6}\nu_{3}$ by $p^{6}\nu_{-3}$ in our initial expression for
$\nu_{0}$ gives%
\[
\tilde{E}q_{4}:\nu_{0}=2p^{6}\nu_{-3}+6p^{5}q\nu_{-2}+15p^{4}q^{2}\nu
_{-1}+20p^{3}q^{3}+15p^{2}q^{4}\nu_{1}+6pq^{5}\nu_{2}.
\]
Also, replacing $q^{2}\nu_{1}$ by $p^{2}\nu_{-1}$ and $q^{4}\nu_{2}$ by
$p^{4}\nu_{-2}$ throughout $\tilde{E}q_{1}$, $\tilde{E}q_{2}$, $\tilde{E}%
q_{3}$ and $\tilde{E}q_{4}$ reduces the number of variables to four. \ The
simultaneous solution is%

\begin{align*}
\nu_{0} &  =\tfrac{-2p\left(  -3-18p+92p^{2}-120p^{3}+816p^{4}-2816p^{5}%
+3840p^{6}-2304p^{7}+512p^{8}\right)  }{{}}\\
&  \tfrac{+2(q-p)^{3}\left(  1+4p+12p^{2}-32p^{3}+16p^{4}\right)
\theta-(q-p)^{2}\theta\sqrt{2}}{{}}\\
&  \tfrac{\cdot\left[  \left(  -1-8p+8p^{2}\right)  \left(  -1-2p-46p^{2}%
-32p^{3}+848p^{4}-2432p^{5}+3200p^{6}-2048p^{7}+512p^{8}\right)  \right.  }%
{{}}\\
&  \tfrac{\left.  +\left(  1+8p+88p^{2}-448p^{3}+1888p^{4}-4864p^{5}%
+6400p^{6}-4096p^{7}+1024p^{8}\right)  \theta\right]  ^{1/2}}{1-256p^{3}%
+768p^{4}-768p^{5}+256p^{6}},
\end{align*}%
\begin{align*}
\nu_{-1} &  =\tfrac{-6pq-(q-p)^{2}\left(  1+32p^{2}-64p^{3}+32p^{4}\right)
\theta+(q-p)\theta}{{}}\\
&  \tfrac{\cdot\left[  \left(  q-p\right)  ^{2}\left(  -1-8p+8p^{2}\right)
\left(  1-8p-56p^{2}-128p^{3}+704p^{4}-768p^{5}+256p^{6}\right)  \right.  }%
{{}}\\
&  \tfrac{\left.  +2\left(  1-32p^{2}+64p^{3}+992p^{4}-4096p^{5}%
+6144p^{6}-4096p^{7}+1024p^{8}\right)  \theta\right]  ^{1/2}}{2p^{2}\left(
1-256p^{3}+768p^{4}-768p^{5}+256p^{6}\right)  },
\end{align*}%
\begin{align*}
\nu_{-2} &  =\tfrac{-6p^{2}q^{2}\left(  5-256p^{3}+768p^{4}-768p^{5}%
+256p^{6}\right)  -(q-p)^{4}\left(  1+8p+24p^{2}-64p^{3}+32p^{4}\right)
\theta+(q-p)\theta}{{}}\\
&  \tfrac{\cdot\left[  \left(  -1+8p+8p^{2}+288p^{3}-944p^{4}-3136p^{5}%
+3776p^{6}-73728p^{7}+712704p^{8}-2445312p^{9}\right.  \right.  }{{}}\\
&  \tfrac{\left.  +4345856p^{10}-4489216p^{11}+2736128p^{12}-917504p^{13}%
+131072p^{14}\right)  }{{}}\\
&  \tfrac{\left.  +2\left(  1+4p+12p^{2}-32p^{3}+16p^{4}\right)  \left(
1-8p+40p^{2}-320p^{3}+1824p^{4}-4864p^{5}+6400p^{6}-4096p^{7}+1024p^{8}%
\right)  \theta\right]  ^{1/2}}{2p^{4}\left(  1-256p^{3}+768p^{4}%
-768p^{5}+256p^{6}\right)  },
\end{align*}%
\begin{align*}
\nu_{-3} &  =\tfrac{-2p\left(  -3-18p-33p^{2}+255p^{3}+441p^{4}-643p^{5}%
-8448p^{6}+28416p^{7}-40448p^{8}+30720p^{9}-12288p^{10}+2048p^{11}\right)
}{{}}\\
&  \tfrac{-(q-p)^{2}\left(  -2-10p-41p^{2}+190p^{3}+89p^{4}-1088p^{5}%
+1728p^{6}-1152p^{7}+288p^{8}\right)  \theta-(q-p)\theta}{{}}\\
&  \tfrac{\cdot\left[  \left(  2+24p+186p^{2}+272p^{3}-1239p^{4}%
-8796p^{5}+43998p^{6}-51924p^{7}-581577p^{8}+3019604p^{9}-4358340p^{10}%
-6991872p^{11}\right.  \right.  }{{}}\\
&  \tfrac{\left.  +38416256p^{12}-72366336p^{13}+79163136p^{14}-54779904p^{15}%
+23804928p^{16}-5971968p^{17}+663552p^{18}\right)  }{{}}\\
&  \tfrac{+2\left(  1+10p+69p^{2}-18p^{3}-330p^{4}-3840p^{5}+8160p^{6}%
+62940p^{7}-250095p^{8}+56320p^{9}+1459360p^{10}\right.  }{2p^{6}\left(
1-256p^{3}+768p^{4}-768p^{5}+256p^{6}\right)  }\\
&  \tfrac{\left.  \left.  -4044096p^{11}+5577696p^{12}-4608000p^{13}%
+2322432p^{14}-663552p^{15}+82944p^{16}\right)  \theta\right]  ^{1/2}}%
{2p^{6}\left(  1-256p^{3}+768p^{4}-768p^{5}+256p^{6}\right)  },
\end{align*}
yielding%
\begin{align*}
\nu_{1} &  =\tfrac{-6pq-(q-p)^{2}\left(  1+32p^{2}-64p^{3}+32p^{4}\right)
\theta+(q-p)\theta}{{}}\\
&  \tfrac{\cdot\left[  \left(  q-p\right)  ^{2}\left(  -1-8p+8p^{2}\right)
\left(  1-8p-56p^{2}-128p^{3}+704p^{4}-768p^{5}+256p^{6}\right)  \right.  }%
{{}}\\
&  \tfrac{\left.  +2\left(  1-32p^{2}+64p^{3}+992p^{4}-4096p^{5}%
+6144p^{6}-4096p^{7}+1024p^{8}\right)  \theta\right]  ^{1/2}}{2\left(
1-256p^{3}+768p^{4}-768p^{5}+256p^{6}\right)  q^{2}},
\end{align*}%
\begin{align*}
\nu_{2} &  =\tfrac{-6p^{2}q^{2}\left(  5-256p^{3}+768p^{4}-768p^{5}%
+256p^{6}\right)  -(q-p)^{4}\left(  1+8p+24p^{2}-64p^{3}+32p^{4}\right)
\theta+(q-p)\theta}{{}}\\
&  \tfrac{\cdot\left[  \left(  -1+8p+8p^{2}+288p^{3}-944p^{4}-3136p^{5}%
+3776p^{6}-73728p^{7}+712704p^{8}-2445312p^{9}\right.  \right.  }{{}}\\
&  \tfrac{\left.  +4345856p^{10}-4489216p^{11}+2736128p^{12}-917504p^{13}%
+131072p^{14}\right)  }{{}}\\
&  \tfrac{\left.  +2\left(  1+4p+12p^{2}-32p^{3}+16p^{4}\right)  \left(
1-8p+40p^{2}-320p^{3}+1824p^{4}-4864p^{5}+6400p^{6}-4096p^{7}+1024p^{8}%
\right)  \theta\right]  ^{1/2}}{2\left(  1-256p^{3}+768p^{4}-768p^{5}%
+256p^{6}\right)  q^{4}},
\end{align*}
in particular. \ 

We must take $\Omega=\{0,-1,-2\}$ as the absorbing set. \ The rate $\lambda$
of clumps of visits to $\Omega$ is equal to $\lambda_{0}+\lambda_{-1}%
+\lambda_{-2}$ where parameters $\lambda_{0}$, $\lambda_{-1}$ and
$\lambda_{-2}$ are solutions of the system%
\[
\lambda_{0}+\lambda_{-1}\nu_{-1}+\lambda_{-2}\nu_{-2}=\left(  1-\nu
_{0}\right)  \pi_{j},
\]%
\[
\lambda_{0}\nu_{1}+\lambda_{-1}+\lambda_{-2}\nu_{-1}=\left(  1-\nu_{0}\right)
\pi_{j+1}\sim\frac{p^{2}}{q^{2}}\left(  1-\nu_{0}\right)  \pi_{j},
\]%
\[
\lambda_{0}\nu_{2}+\lambda_{-1}\nu_{1}+\lambda_{-2}=\left(  1-\nu_{0}\right)
\pi_{j+2}\sim\frac{p^{4}}{q^{4}}\left(  1-\nu_{0}\right)  \pi_{j}.
\]
The total clump rate is consequently%
\[
\lambda\sim\frac{(q-p)u+(q-p)^{3}\theta+\sqrt{2}(q-p)^{2}\sqrt{v+u\theta}%
}{4q^{4}}\pi_{j}%
\]
where $u$ and $v$ appear at the end of Section 3.1. \ Thus the exponential
coefficient is%
\[
\varepsilon_{1}=\frac{1}{6}\cdot\frac{p^{2}}{q^{2}}\cdot\frac{q^{3}}{p^{3}%
}\cdot\frac{\lambda}{\pi_{j}}\cdot A(p)\sim\frac{\left[  (q-p)u+(q-p)^{3}%
\theta+\sqrt{2}(q-p)^{2}\sqrt{v+u\theta}\right]  ^{2}}{288pq^{9}}=\frac
{\chi_{3}(p)}{6}%
\]
where%
\[
\chi_{3}(p)=\frac{(q-p)^{2}}{12pq^{9}}\left[  \alpha+(q-p)^{2}\beta
\theta+(q-p)\sqrt{2}\sqrt{\gamma+\alpha\beta\theta}\right]  ,
\]%
\[
\alpha=1-4p+10p^{2}-52p^{3}+226p^{4}-520p^{5}+640p^{6}-400p^{7}+100p^{8},
\]%
\[
\beta=1-2p+6p^{2}-8p^{3}+4p^{4},
\]%
\begin{align*}
\gamma &  =1-4p+16p^{2}-104p^{3}+506p^{4}-1808p^{5}+5604p^{6}-15576p^{7}%
+35574p^{8}\\
&  -61160p^{9}+75152p^{10}-63440p^{11}+34840p^{12}-11200p^{13}+1600p^{14}%
\end{align*}
was conjectured in \cite{Fi1-heu}. \ Again, what was missed previously is the
expression for $\varepsilon_{0}$ as a perfect square, a corollary of the
hidden relation
\[
\frac{\lambda}{\pi_{j}}=3q^{2}A(p).
\]
We conjecture, for arbitrary $\ell\geq2$, that%
\[
\frac{\lambda}{\pi_{j}}=\ell q^{2}A(p)
\]
among many possible unsolved problems.\footnote[5]{In this paper, we imagined
a single line of traffic flowing from east to west. \ In the following,
imagine instead two independent Bernoulli($p$) queues with no $\ell$ parameter
at all. Think of one queue for east-to-west traffic (EW) and one queue for
north-to-south (NS) traffic. \ The intersection of the two lines of traffic is
governed by one stoplight. \ An (old) green light for EW traffic turns red
precisely when then are no EW cars left, the (new) green light for NS traffic
turns red precisely when then are no NS cars left, and so on. \ What is the
stationary distribution for this scenario? \ This may be a difficult question
to answer. \ From \cite{Ald-heu}, it seems that, if $\pi$ can be calculated,
then the distribution of the maximum line length of idle EW cars (say) can be
readily found via the clumping heuristic.}\footnote[6]{Another
conceptualization of the $\ell=2$ scenario involves patients randomly arriving
at a hospital emergency room. \ One doctor treats a new patient starting at
times $\equiv3\operatorname{mod}4$; the other likewise at times $\equiv
4\operatorname{mod}4$. \ Services lengths are constant. \ Allowing these to be
nondeterministic complicates the analysis.}\bigskip

\textbf{Acknowledgement} Our reliance in Sections 2 and\ 3 on the computer
algebra systems Mathematica and Maple is obvious; we thank all people who
created such software.%

\begingroup\renewcommand{\enotesize}{\normalsize}
\setlength{\parskip}{2ex}
\theendnotes\endgroup

\end{document}